\documentclass[12pt]{article}
\usepackage{amsthm,amsmath,amssymb,latexsym,graphics}
\begin{document}
\newtheorem{lem}{Lemma}[section]
\newtheorem{rem}{Remark}[section]
\newtheorem{question}{Question}[section]
\newtheorem{prop}{Proposition}[section]
\newtheorem{cor}{Corollary}[section]
\newtheorem{thm}{Theorem}[section]
\newtheorem{problem}{Problem}[section]

\title
{A geometric problem and the Hopf Lemma. I
}
\author{YanYan Li\thanks{Partially
 supported by
      NSF grant DMS-0401118.}
\\
			            Department of Mathematics\\
				           Rutgers University\\
					          110 Frelinghuysen Road\\
						         Piscataway, NJ 08854\\
							        USA
\\
\\
Louis Nirenberg\thanks{With thanks to the
Pacific Institute of Mathematical Sciences for support.}
 \\
Courant Institute\\
251 Mercer Street \\
New York, NY 10012\\
\\
\\
\\
Dedicated to Antonio Ambrosetti on his 60th birthday
\thanks{ This is a report of a talk given at the January
2005 conference in Rome in celebration of his birthday}
}

\date{ }
\maketitle

\input { amssym.def}

\setcounter {section} {0}

\section{Introduction}

The problem we consider starts with the following classical result of A.D.
Alexandrov \cite{A}:

\begin{thm} (\cite{A})\ Let $M$ be a compact smooth hypersurface,
without boundary, embedded in $\Bbb R^{n+1}$
with the property that the mean curvature (average
of principle curvatures, using interior normal)
is identically constant.  Then $M$ is a sphere.
\label{thma} 
\end{thm}

If $M$ is immersed instead of embedded, the conclusion
of the theorem may fail, even in dimension
$n=2$.  Indeed, 
in 1986, Wente \cite{W} constructed a counter example
in case $M$ is an immersed torus, with
self intersection, in $\Bbb R^3$.
A. Ros \cite{Ros} in 1987 extended Theorem \ref{thma}
from mean curvature to the elementary
symmetric functions of the principal
curvatures of $M$.
In 1997 YanYan Li \cite{L} gave some far reaching generalizations  
including very general symmetric functions of the principal curvatures of 
$M$.  But here we just mention one of the results --- still for
the  mean curvature.

\begin{thm} (\cite{L})\ Let $M$ be a compact smooth hypersurface
without boundary embedded in $\Bbb R^{n+1}$.  Let $K$
be
a $C^1$ function in $\Bbb R^{n+1}$ satisfying
\begin{equation}
\frac{\partial K}{\partial x_{n+1}}\le 0,
\label{1}
\end{equation}
Suppose that at each point $x$ of $M$ the mean curvature $H(x)=K(x)$.
Then 
$M$ is symmetric about some hyperplane
\begin{equation}
x_{n+1}=\lambda_0.
\label{2}
\end{equation}
\end{thm}

Li then proposed that we consider the more 
general question in which the condition
$H(x)=K(x)$ with $K$ satisfying 
(\ref{1}),
is replaced by the weaker, more natural,
condition: Whenever $(x', a)$ and $(x', b)$, $a<b$, lie on $M$
(here $x'=(x_1, \cdots, x_{n})$) then
\begin{equation}
H(x', b)\le H(x', a).
\label{3}
\end{equation}

\begin{question} Is it true then that $M$ is symmetric about some hyperplane
$x_{n+1}=\lambda_0$?
\label{question1}
\end{question}

This paper --- here we consider only 
one-dimensional problems --- and
its sequel, are concerned with
this question.

First we recall Alexandrov$'$s argument.
It introduces the, now familiar, method of moving planes,
and the proof replies on the strong maximum principle and
the Hopf Lemma for second order elliptic
equations.  Here it is:  $M$ is the boundary of an open set
$U$ in $\Bbb R^{n+1}$. For $\lambda$ less than, but close to,
$
\displaystyle{
\max_M x_{n+1},
}
$
take the part $S_\lambda$
of $M$ lying above (i.e. larger $x_{n+1}$)
$\lambda$ and reflect it in the plane
$
x_{n+1}=\lambda.
$
The reflected piece of 
surface, $S_\lambda'$, lies in 
$\overline U$.  Decrease $\lambda$ and continue
to reflect $S_\lambda$ so that
$S_\lambda'$ continues to lie in $\overline U$.  There
will be a first value $\lambda_0$ of $\lambda$, such that one
of the two things happen:

(i)$S_{\lambda_0}'$ touches $M$ at some point
$(x_0', a_0)$ with $a_0<\lambda$, and the line
$\{(x_0', x_{n+1})\ |\
x_{n+1}\in \Bbb R^{n+1}\}$
is transversal to $M$ at $(x_0', a_0)$,

or

(ii)\  At some point $P$ on $x_{n+1}
=\lambda_0$ the hypersurfaces $S_{\lambda_0}$ and $M$ are 
tangent to each other.

\medskip

Note that both
 things may happen at the same
$\lambda_0$.

In Case (i) we may describe $M$ and $S_{\lambda_0}'$
near $(x_0', a_0)$
as graphs of smooth functions $v(x')$, $u(x')$ with 
\begin{equation}
v(x')\le u(x')
\ \ \mbox{and}\ v(x_0')=u(x_0').
\label{4}
\end{equation}
Both functions satisfy, near $x_0'$, the nonlinear elliptic
equation expressing the fact that the mean 
curvature is the constant $H$, 
$$
H[u]:= \nabla \left( \frac{\nabla u}{  \sqrt{1+|\nabla u|^2 }  }\right)=H.
$$
But by the strong maximum principle it follows
from 
(\ref{4}) that $u\equiv v$ near $x_0'$.  Applying 
this argument at other points on $M$ it follows that
$$
S_{\lambda_0}'=\{x\in M\ |\ x_{n+1}<\lambda_0\}.
$$
This is the desired symmetry.

In Case (ii) we turn the picture around

***********

Fig. 1

***********

In these coordinates, call them $y$, with
$y_1=P_1-x_{n+1}$, 
we may represent $S_{\lambda_0}'$, near $P$, as the
graph of a function
$u(y')$, and the part of $M$ lying in 
$y_1>0$ as a graph of some $v(y')$.  We have, with some abuse of notation, 
$$
u(y')\ge v(y'), \ \ 
u(P)=v(P),\ \ \mbox{and}\
\nabla u(P)=\nabla v(P).
$$
As before,
$u$ and $v$ satisfy the same elliptic equation
in $y_1>0$, near the origin.  By the Hopf
Lemma, $u\equiv v$ near the origin.  Then,
using the strong maximum 
principle we extend this fact globally, to conclude the desired symmetry.

In \cite{L} Li uses the moving plane method but makes essential use of the fact that the function $K$ is locally Lipschitz in $\Bbb R^{n+1}$.

What happens if, following Alexandrov,
we try to use moving planes for the problem where
$H(x', b)\le H(x', a)$ for $b>a$?
We are led again to the two cases (i) and (ii)
above.  Case (i) is easily handled.  Here is the picture

***********

Fig.2 

***********

Again we have two functions $u(x')\ge v(x')$.
But now $H[u]\le H[v]$.  We may
still use the strong maximum principle
and infer that $u\equiv v$.

The trouble arises
in Case (ii).  If we redraw
Fig. 1,
we have

***********

Fig. 3 

***********

We may represent $S_{\lambda_0}'$ and $M$, in $y_1>0$ by $u(y')$ and
$v(y')$, with $u\ge v$.  However the condition that
the mean curvature of $S_{\lambda_0}'$ at $B$ is
$\le $ that at $A$,  compares the mean curvature
of $u$ and $v$ $\underline{  \mbox{but at
different points} }$,
$(y_1, y'')$ and $(\bar y_1, y'')$ where
$u(y_1, y'')= v(\bar y_1, y'')$,
with $y_1\le \bar y_1$.

Thus we are led to look for a more general 
form of the Hopf Lemma.

Before stating some results, we point out that
the answer to Question \ref{question1} is
false
in general, even for a closed curve with
interior convex in the $x_2-$direction --- as 
 the following figures shows.

***********

Fig. 4 

***********

Here the ends are symmetric to each other, the bottom
bump is symmetric, as is the top bump, 

Even if $M$ is not necessarily symmetric,
would the inequality on the  mean curvature implies that equality holds
in the following sense?

\begin{question}
Is it true that (\ref{3}) implies
that for any
 $A, B\in M, A_{n+1}<B_{n+1}$,
we must have $H(A)=H(B)$?
\label{question2}
\end{question}

In Section 6 we give a counterexample.
However we do not know the answer to
\begin{question}
Is the answer to Question \ref{question2} yes in case we consider
 $A, B\in M$ such that
 for all $0<t<1$,
  $tA+(1-t)B$ lies  inside $M$ ?
 \label{question3}
 \end{question}

In Part II, \cite{LN},  we present our results on 
Question \ref{question1}.  We assume that the 
(embedded) hypersurface is smooth and satisfies

\noindent{\bf Condition S.}\ $M$ stays
on one side of any hyperplane parallel to the $x_{n+1}-$axis that is
tangent to $M$.

We
 believe that this
should suffice to prove symmetry.  However
our proof requires a further condition:

\noindent {\bf Condition T.}\
Any line parallel to the $x_{n+1}-$axis 
that is tangent to $M$ has contact of finite order.

Condition T automatically holds in case
$M$ is analytic; while Condition S automatically holds in case
$M$ is convex.

Our main result is
\begin{thm} (\cite{LN})\  Let $M$  be a smooth  compact embedded  hypersurface
in $\Bbb R^{n+1}$ satisfying 
$$
H(x', x_{n+1})\le H(x', \tilde x_{n+1})
$$
 for any two points $ (x', x_{n+1}), (x', \tilde x_{n+1})\in M$
 satisfying $x_{n+1}\ge \tilde x_{n+1}$.
 Then, if Condition S and Condition T hold,
  $M$ must be symmetric with respect to some hyperplane
   $x_{n+1}=constant$.
  \label{thm1.3}
 \end{thm}

In this paper we restrict ourselves
to curves. 
The main result in Part I is
\begin{thm} Let $M$ be a closed $C^2$  embedded curve 
in the plane  satisfying Condition S.  Assume that whenever $(x_1, a)$,
$(x_1, b)$, with $a<b$,
lie on the curve,
\begin{equation}
\mbox{curvature of}\ M\
\mbox{at}\
(x_1, b)\le 
\mbox{curvature of}\ M\
\mbox{at}\
(x_1, a).
\label{2.1}
\end{equation}
Then $M$ is symmetric about some line
$x_2=\lambda_0$.
\label{thm1}
\end{thm}

\begin{rem} In Theorem \ref{thm1}, we do not assume Condition T.
\end{rem}

The theorem is proved in Section 3.  In Part II,
in addition to mean curvature, we also extend
Theorem \ref{thm1.3} to other symmetric functions of the principal
curvatures.  A number of open problems
are also presented there.  In Section 7
of this part I we also mention several
which are local in nature.

\section{One dimensional model problems}

We first looked at some very simple one-dimensional model
problems which seemed to us to be of interest.
Here is one of them.

\begin{thm}
Let $u\ge v$ be positive $C^2$ functions on 
$(0,b)$, which are also
in $ C^1( [0, b])$.
  Assume that
$$
u(0)=\dot u(0)=0,
$$
and
\begin{equation}
\mbox{either}\  \dot u>0\ \mbox{on}
\ (0, b)\ \mbox{or}\
\ \dot v>0 \ \mbox{on}\ (0,b).
\label{eitheror}
\end{equation}
Our main hypothesis is :
\begin{equation}
\mbox{whenever }\ u(t)=v(s)\ \mbox{for}\
0<t<s<b\ \mbox{we have}
\ \ddot u(t)\le v''(s).
\label{2.2}
\end{equation}
(here $\cdot = \frac d{dt}$, $'=\frac{d}{ds}$).
Conclusion:
$$
u\equiv v\ \ \mbox{on}\ [0, b].
$$
\label{prop1}
\end{thm}

\begin{rem}
This is a kind of extension of the Hopf
Lemma, for if in place of (\ref{2.2}) we assumed
\begin{equation}
\ddot u(t)\le \ddot v(t)\ \mbox{on}\ (0,b)
\label{2.3}
\end{equation}
the result would simply follow from the Hopf Lemma.
\end{rem}

\begin{rem} If we replace (\ref{eitheror}) by
both $\dot u\ge 0$ and $\dot v\ge 0$ on $(0, b)$, the conclusion of
Theorem \ref{prop1} may fail.  See the following example.
\label{rem1-4}
\end{rem}

\noindent{\bf Example 1.1}\
Let $u\in C^\infty([0, 2])$ satisfy
$$
u(t)=
\left\{
\begin{array}{ll}
\displaystyle{
t^3,}
&
\displaystyle{
0\le t\le \frac 13
}\\
\displaystyle{
1+(t-1)^3,
}
&
\displaystyle{
\frac 23\le t\le 2,
}
\end{array}
\right.
$$
$$
\dot u>0\quad \mbox{in}\ (0,1),
$$
and let
$$
v(t)=w(t)=
\left\{
\begin{array}{rl}u(t), & 0\le t\le 1,\\
1,& 1< t\le 2,
\end{array}
\right.
$$
as in

***********

Fig. 5 

***********

Before proving Theorem \ref{prop1}, we give a few lemmas.
Some of these are not really used in the proof of Theorem \ref{prop1},
but seem of interest.

First, a variation of the strong maximum principle.
\begin{lem} Let $u\ge v$ be two $C^2$ functions on
$(0,b)$.  We assume (\ref{2.2}) and
\begin{equation}
\max\{\dot u, \dot v\}>0\ \mbox{on}\ (0,b).
\label{aa0}
\end{equation}
Then either
$$
u>v\  \mbox{on}\ (0,b),
$$
or
\begin{equation}
u\equiv v\ \ \mbox{on}\ (0, b).
\label{aa1}
\end{equation}
\label{lema1}
\end{lem}

\begin{rem} If we replace (\ref{aa0})
by $\dot u, \dot v\ge 0$ on
$(0, b)$, the conclusion of the lemma
may fail.
See Example 1.1 above.
\label{rem1-1}
\end{rem}

\noindent{\bf Proof of Lemma \ref{lema1}.}\
Suppose $u(c)=v(c)$ for some $0<c<b$. Then, by (\ref{aa0}),
$$
\dot u(c)=\dot v(c)>0.
$$
By  the
implicit function theorem, for $s$ close to $c$, there is a $C^2$ function $t(s)$ such that
\begin{equation}
u(t(s))=v(s).
\label{2.4}
\end{equation}
For $s$ close to $c$, set
$$
g(s)=s-t(s),
$$
so $g\ge 0$.  Differentiating
(\ref{2.4}) we find
$$
\dot u(t(s)) t'=v'(s),
\quad \ddot u t'^2 +\dot u t''=v''.
$$
In terms of $g$ the last equation becomes,
by (\ref{2.2}),
\begin{equation}
0\ge \ddot u(t(s))-v''(s)
=\dot u g'' -g'(g'-2)\ddot u.
\label{2.5}
\end{equation}

Now if $u(s)=v(s)$ for some $s>0$
then $g$ vanishes there.  By the strong maximum
principle applied to (\ref{2.5})
 implies that $v\equiv u$ in a neighborhood of $c$.
 By the same argument, $u\equiv v$ in
 a larger neighborhood, and 
(\ref{aa1}) then follows.
Lemma \ref{lema1} is established.

\vskip 5pt
\hfill $\Box$
\vskip 5pt

\begin{lem}
Let $u$ and $w$ be positive $C^{1,1}_{loc}$ functions 
on $(0,c)$, belonging to $C^1( [0, c))$, and satisfying,
for some $f\in L^\infty_{loc}(0,\infty)$,
$$
\ddot u=f(u), \ \ \ 
\ddot w=f(w),\qquad \mbox{in}\ (0,c),
$$
and
$$
u(0)=w(0)=0,
\ \ \
\dot u(0)=\dot w(0),\ \ \
\dot u>0\ \mbox{on}\ (0,c).
$$
Then
$$
u\equiv w\ \mbox{on}\ (0,c).
$$
\label{lem1prime}
\end{lem}

\begin{rem}  We do not assume $\dot w>0$.
On the other hand, if we replace $\dot u>0$
by $\dot u\ge 0$ in the hypotheses,  
the conclusion may fail.
See Example 1.1.
\label{rem1-2}
\end{rem}

Here is another  simple uniqueness result; it could be taught
in a beginning course on ordinary differential equations.

\begin{lem} Let $u$ and $w$ be positive $C^{1,1}_{loc}$ functions
in $(0, c)$, belonging to $C^1([0, c))$, and both satisfying
\begin{equation}
\ddot u=f(u),
\label{2.8}
\end{equation}
and
$$
u(0)=w(0)=0,
\quad
\dot u(0)=\dot w(0).$$
Assume that $f(\rho)$ is locally Lipschitz for
$\rho>0$ (not necessarily on $\rho\ge 0$).
Then
$$
u\equiv w.
$$
\label{lem1}
\end{lem}

Note that we assume neither  $\dot u\ge 0$ nor $\dot w\ge 0$.

\noindent{\bf Proof
of Lemma \ref{lem1}.}\ 
The proof is by obtaining an ``explicit'' expression for $u(t)$.
Multiply (\ref{2.8}) by $2\dot u(t)$ and integrate
from some $t_0>0$ to some $t>t_0$.  We
find
$$
\dot u(t)^2 -\dot u(t_0)^2=F(u(t))-
F(u(t_0)).
$$
Here $F(\rho)$ is such that
$$
\frac{dF}{d\rho}=2f(\rho)
\qquad\mbox{for}\ \rho>0.
$$
Letting $t_0\to 0$ we see that $F(\rho)$ has a limit as $\rho\to 0$, which
we may take to be $\dot u(0)^2$.  Thus, letting
$t_0\to 0$ we find
\begin{equation}
\dot u(t)^2 =F(u(t)).
\label{2.13}
\end{equation}

\noindent $\underline{  \mbox{Claim} }.$\ On $(0, c/2)$, $\dot u>0$.

For if not, if $\dot u(t_1)=0$ for some
$0<t_1<c/2$, we would have, by the local
Lipschitz property of $f$ for $\rho>0$, that
the function $u$ is symmetric about $t_1$.
But then it would have to vanish at
$2t_1$ --- where $u$
is positive.

Consequently, from (\ref{2.13}), we find that
$F(\rho)>0$ for $0<\rho< u(\frac c4)$, and
$$
\dot u(t)=\sqrt{  F(u(t))  }
$$
or 
\begin{equation}
\frac{  \dot u }{  \sqrt{F(u)}  }=1.
\label{2.14}
\end{equation}
If, on $\rho>0$, $G(\rho)$ is such that
$$
G_\rho=\frac 1{   \sqrt{F(\rho)}  },
$$
we find from (\ref{2.14}) that
$$
\frac {d}{dt} G(u)=1.
$$
Integrating from $t_0>0$ to $t>t_0$ we obtain
$$
G(u(t))-G(u(t_0))=t-t_0.
$$
Letting again $t_0\to 0$
we see that $G$ has a limit at $\rho=0$, which we may take to be
zero.  Thus
$$
G(u(t))=t.
$$
So $u(t)$ is uniquely determined on $(0, c/2)$,
since $G_\rho>0$ for $\rho>0$.
Then, by the local Lipschitz property of $f(\rho)$ on
$\rho>0$, it follows that $u$ is unique on 
$(0, c)$.  Lemma
\ref{lem1} is proved.

\vskip 5pt
\hfill $\Box$
\vskip 5pt

 \noindent{\bf Proof of Lemma \ref{lem1prime}.}\
Follow the proof of Lemma \ref{lem1} until (\ref{2.13}). 
Similarly we also have
\begin{equation}
\dot w(t)^2 =F(w(t))\qquad \mbox{on}\ (0,c).
\label{abc1}
\end{equation}
Let $0<b<c$ be any number satisfying
\begin{equation}
\max_{ [0, b] }w<\sup_{ [0, c) }u
=\lim_{t\to c^-}u(t).
\label{number}
\end{equation}
We will prove that
\begin{equation}
u\equiv w\quad \mbox{on}\ [0, b].
\label{A4-1}
\end{equation}

By  (\ref{2.13}), (\ref{abc1}) and the fact that $\dot u>0$ on $(0, c)$, 
we know that
$$
\dot w(t)^2=F(w(t))>0,
\qquad\forall\ 0<t<b.
$$
Since $w(0)=0$ and $w>0$
on $(0, b)$, we have, in view of the above, $\dot w>0$ on $(0, b)$.
Proceed as in the proof of Lemma \ref{lem1}, we arrive at
$$
G(u(t))=G(w(t))=t,\qquad
0<t<b.
$$
But $G_\rho>0$ and $G(0)=0$, we obtain
as before (\ref{A4-1}).

Arguing in the same way we see that 
$u\equiv w$ on an interval $(0, b')$, $b'>b$ ---
and then $(0, c)$.

\vskip 5pt
\hfill $\Box$
\vskip 5pt

The following lemma can be viewed as a variation of the
maximum principle.
\begin{lem} Let $u, v\in C^2((0,b))\cap C^0([0, b])$
satisfy
$$
v(0)\le u(0),\
\mbox{and}\ v< u(b)\ \mbox{on}\ (0,b),
$$
\begin{equation}
\mbox{either}\ \dot u>0\  \mbox{on}\ (0, b),\
\mbox{or}\ \dot v>0 \ \mbox{whenever}\ u(0)<v<u(b),
\label{either}
\end{equation}
and
(\ref{2.2}).
Then
$$
u\ge v\ \mbox{on}\ [0,b].
$$
\label{maximum}
\end{lem}

\begin{rem} If we change (\ref{either}) to $\dot u\ge 0$
and $\dot v\ge 0$, the
conclusion of Lemma \ref{maximum} may fail.
See the example below.
\label{rem1-3}
\end{rem}

\noindent{\bf Example 1.2.}\
Let $u\in C^\infty([0,4])$ satisfy
$$
u(t)=
\left\{
\begin{array}{ll}1, & 0\le t\le 2,\\
\displaystyle{
1+(t-2)^3,
}
&
\displaystyle{
2< t\le \frac 73,}\\
2, & 
\displaystyle{
\frac 83< t\le 3,}\\
\displaystyle{
2+(t-3)^3,}
& 3<t\le 4,
\end{array}
\right.
$$
$$
\dot u>0\quad \mbox{in}\ (\frac 73, \frac 83),
$$
and let $w\in C^\infty([1, 5])$,
with $w'$ nonnegative, satisfy
$w(1)=0$,
$$
w(t)=
\left\{
\begin{array}{rl}u(t), & 2\le t\le 3,\\
2,& 3<t\le 5.
\end{array}
\right.
$$
Then let
$$
v(t)=w(t+1),\qquad 0\le t\le 4.
$$
See

***********

Fig. 6 

***********

\noindent{\bf Proof of Lemma \ref{maximum}.}\
Shift $v$ far to the right and then slide it to the left
until its graph first touches that of $u$.
If the touching occurs at $(c, u(c))$ for some $0<c<b$,
then $u'(c)>0$, and therefore, by  Lemma \ref{lema1}, 
$u$ and the shift of $v$ must coincide
near $c$.
Again by the lemma, the
set of points where $u\equiv v$ is open.  Since it is
closed, we conclude that the shift of $v$ is
$v$ itself and $u\equiv v$.  
Otherwise, we can slide the shift of $v$ all the way to the origin and
we conclude $u\ge v$ on $(0, b)$.

\vskip 5pt
\hfill $\Box$
\vskip 5pt

\noindent{\bf Proof of Theorem
\ref{prop1}.}\ (a)\ We first assume that
$\dot u>0$ on $(0,b)$.
Because of  Lemma \ref{lema1},  
 we may suppose that 
$$
u(t)>v(t)\quad
\mbox{for}\ t>0,
$$
and we will derive a contradiction.

Our proof makes use of the fact that $u$
satisfies some differential equation.  Namely,
since $\dot u(t)>0$ for $t>0$, we see that for
$u>0$, $t$ is a $C^2$ function of $u$.
It follows that we may write
\begin{equation}
\ddot u=f(u),
\label{2.6}
\end{equation}
with $f$ some unknown function which is however  
continuous in $u$ on $[0, u(b)]$.
The main hypothesis
(\ref{2.2}) is then equivalent to
the following
for $v$:
\begin{equation}
v''(s)\ge f(v(s)).
\label{2.7}
\end{equation}
Thus $v$ is a subsolution of (\ref{2.6})
while $u\ge v$ is a solution.

Consequently there is a solution
$w$ of (\ref{2.6}) lying between $u$ and $v$, with
$$
w(0)=0, \  \
w(c)=v(c)\quad
\mbox{for some fixed}\ c\
\mbox{in}\ (0,\bar b).
$$
If $\ddot u$ is locally Lipschitz on $(0,b)$, then 
 $f$ is locally Lipschitz continuous in $(0, u(b)]$ and  we
 can see this by first, for small positive
$\epsilon$, finding a solution $w_\epsilon$ between $u$ and $v$
on $[\epsilon, c]$ with
$$
w(\epsilon)=v(\epsilon), \ \ w(c)=v(c).
$$
Because of  the locally  Lipschitz property of $f$,  
the usual argument by monotone
iterations yields $w_\epsilon$.  Letting $\epsilon\to 0$
one easily obtains a solution $w$ satisfying
(\ref{2.6}) with
$$
v(t)\le w(t)\le u(t).
$$
J. Mawhin pointed out to us that the
approximation by $w_\epsilon$ is not necessary.  That there is
a solution between $u$ and $v$, even for merely continuous
$f$, is known, see
\cite{M}.

We now have two positive solutions in $(0, c)$, $u$ and
$w$ of (\ref{2.6})
with $u(0)=\dot u(0)=w(0)=\dot w(0)=0$.
By Lemma \ref{lem1prime}, $u\equiv w$ on $(0, c)$, violating
$u(c)>v(c)$. Impossible.

(b)\ We now assume that $\dot v>0$ on $(0, b)$.
Let $0<a<b$ be any number satisfying
\begin{equation}
\max_{ [0, a]}u< \sup_{ [0, b)}v
=\lim_{t\to b^{-}}v(t).
\label{B2-0}
\end{equation}
We will prove that
\begin{equation}
u\equiv v\ \mbox{on}\ [0, a].
\label{B2-1}
\end{equation}
It is easy to see that this would imply
$u\equiv v$ on $[0, b]$.

Since $\dot v>0$ on $(0, b)$ it follows, in view of
(\ref{B2-0}), that for every $t\in (0, a)$ there
is a $C^2$ function $s(t)$ such that
\begin{equation}
u(t)=v(s(t)),\qquad
0<t<a.
\label{B3-1}
\end{equation}
Set
$$
g(t)=s(t)-t,\quad 0<t<a,
$$
so $g\ge 0$.  Differentiating (\ref{B3-1}) we find,
still use notation
 $\cdot = \frac d{dt}$, $'=\frac{d}{ds}$,
 $$
 \dot u(t)=v'(s(t))\dot s(t), \quad
 \ddot u=v'' \dot s^2+v' \ddot s.
 $$
 In terms of $g$ the last equation becomes, by (\ref{2.2}),
 $$
 0\ge \ddot u(t)-v''(s(t))\ge
 v' \ddot g+\dot g(\dot g+2)v''.
 $$
 If $u(t)=v(t)$ for some $0<t<a$, then $g$ vanishes there and,
 as before, $g\equiv 0$ on
 $(0, a)$ which in turn implies (\ref{B2-1}).  

 Thus we may assume that
 $$
 u>v\quad \mbox{on}\ (0, a),
 $$
 and we will derive a contradiction.

 Since $v'>0$ on $(0, b)$, we may, as before, write
 \begin{equation}
 v''=f(v)\quad\mbox{on}\ (0, b)
 \label{B6-1}
 \end{equation}
 where $f$ is some unknown continuous function on
 $[0, \lim_{s\to b^-}v(s))$.
 By our main hypothesis  (\ref{2.2}),
 in view of (\ref{B2-0}),
 $$
 \ddot u\le f(u)\quad \mbox{on}\ (0, a).
 $$
 As before there exists a solution $w$ of
 (\ref{B6-1}) lying between $u$ and $w$, with
 $$
 w(0)=0, \quad w(a)=v(a).
 $$
 By Lemma \ref{lem1prime}, $u\equiv w$ on
 $[0, a]$, violating $u(a)>v(a)$.
Theorem \ref{prop1} is established.

\vskip 5pt
\hfill $\Box$
\vskip 5pt

\section{Proof of Theorem \ref{thm1}}

We first give the main lemma for the proof of Theorem \ref{thm1}.
\begin{lem} Let, for some $b>0$,  $u, v\in C^2( (0, b))\cap C^1([0,b])$
satisfy
$$
u(t)\ge v(t)>0, \qquad \mbox{for}\
0<t\le b,
$$
$$
\mbox{either}\  \dot u(t)>0, \dot v(t)\ge 0\ \mbox{for}\
0<t\le b,\  \mbox{or}\ \dot u(t)\ge 0,  \dot v(t)>0\ \mbox{for}\
0<t\le b,
$$
and
$$
u(0)=\dot u(0)=0.
$$
Assume
\begin{equation}
\mbox{whenever}\ u(t)=v(s)\ \mbox{for}\
0<t<s<b\ \mbox{we have}\
\frac{  \ddot u(t)  }
{ (1+\dot u^2)^{\frac 32}  }
\le \frac{  v''(s) }
{  (1+v'^2)^{ \frac 32}  }.
\label{3.1}
\end{equation}
Then
\begin{equation}
u\equiv v\ \ \mbox{on}\ [0, b].
\label{3.2new}
\end{equation}
\label{lem8}
\end{lem}

\noindent{\bf Proof.}\  We will only prove it under 
``$ \dot u(t)>0, \dot v(t)\ge 0$ for
$0<t\le b$''.  The changes needed when
assuming instead 
``$\dot u(t)\ge 0,  \dot v(t)>0$ for
$0<t\le b$'', are similar to those 
in the proof of Theorem \ref{prop1}.
We start as in the proof of Theorem 
\ref{prop1}.  $u$ satisfies
an equation of the form
\begin{equation}
\frac{  \ddot u}{  (1+\dot u^2)^{ \frac 32}  }=f(u)
\label{3.3}
\end{equation}
with some unknown function $f$ which is
however continuous on $(0,  u(b)]$.

Our condition (\ref{3.1}) means that
\begin{equation}
\frac{  \ddot v  }{  (1+\dot v^2)^{ \frac 32}  }
\ge f(v).
\label{3.4}
\end{equation}

Multiply (\ref{3.3}) by $\dot u$, we find
that 
$$
-
\frac{d}{dt}\left(  \frac 1{  (1+\dot u^2 )^{\frac 12} }\right)
=\frac{d}{dt} F(u)
$$
where $F$ is such that $F_u=f(u)$.
  
  Integrating this from $t_0$ to $t$,
  $t_0>0$, we find 
  $$
  \left(\frac  1{   1+\dot u(t_0)^2)^{\frac 12}  }\right)^{\frac 12}
  -
  \left(\frac  1{   1+\dot u(t)^2)^{\frac 12}  }\right)^{\frac 12}
=F(u(t))-F(u(t_0)).
$$
Letting $t_0\to 0$ we see that
$F(u)$ has a limit at $u=0$, which we may take
to be zero.  Thus
$$
\left(  \frac 1{  1+\dot u^2}\right)^{\frac 12}
=1-F(u)
$$
so that 
\begin{equation}
\dot u=\left[
\frac 1{  (1-F(u))^2 }-1\right]^{\frac 12}.
\label{3.5}
\end{equation}

Next, multiplying (\ref{3.5}) by $\dot v\ge 0$ we obtain
$$
-\frac d{dt}\left(  \frac 1{  1+\dot v^2  }\right)^{ \frac 12}
\ge \frac d{dt} F(v) \ \ 
\mbox{for}\ 0<t_0<t.
$$
Since $F(0)=0$, we find, on integrating,
$$
1-\left(  \frac 1{1+\dot v^2 }\right)^{ \frac 12}\ge F(v(t)).
$$
Thus, since $\dot v\ge 0$,
\begin{equation}
\dot v\ge
\left[  \frac 1{  (1-F(v(t)))^2 }  -1\right]^{\frac 12}.
\label{3.6}
\end{equation}

But (\ref{3.5}) and (\ref{3.6}) imply that
\begin{equation}
\mbox{whenever}\  u(t)=v(s)\
\mbox{for}\ t<s,
\mbox{there} \ \dot u(t)\le v'(s).
\label{3.7}
\end{equation}

Since $\dot u>0$ for $t>0$, by the
implicit function theorem, there is a $C^1$ function
$t\le s$ such that
$$
u(t(s))=v(s).
$$
Thus if $g=s-t\ge 0$, we have
by differentiating,
$$
\dot u(t)( 1-g')=v'(s).
$$
From (\ref{3.7}) it follows that
$$
\dot ug'=\dot u-v'\le 0
$$
i.e. $g'\le 0$.  Since
$g(0)=0$ and $g\ge 0$, it follows that
$g\equiv 0$ --- which implies
(\ref{3.2new}).

\vskip 5pt
\hfill $\Box$
\vskip 5pt

\noindent{\bf Proof of Theorem \ref{thm1}.}\
Condition S implies that there are just two lines parallel
to the $x_2-$axis which are tangent to $M$. 
We carry out
the moving plane method as described in Section 1
$\underline{  \mbox{except that}  }$
we define $\lambda_0$ to be the first value of $\lambda$,
as we decrease it, such that for any
$\lambda<\lambda_0$, $S_\lambda'$
does not lie in $\overline U$.
We then obtain as in Section 1,
cases (i) and (ii).

If Case (i) happens, then it can be  treated
as described there, but with some difference:
flat vertical segments may occur, though we 
still obtain symmetry.  See picture below.

***********

Fig. 7 

***********

Now we look at  Case (ii).
There is a common tangency point
of $S_{\lambda_0}'$ and
$M$ such that if we rotate the figure it looks as follows, with coordinates
$t$ and $y$
and, due to Condition S,
 the curves lie above the
$t-$axis.

***********

Fig. 8 

***********

***********

Fig. 9 

***********

Again, we are allowing $M$ to have some flat segment.  Take as origin, the point $P$
as shown. 
Let $[0, a]$ be the largest interval
for the flat segment of $S_{\lambda_0}'$ and let
$[0, b]$ be the largest interval
for the flat segment of $M$ as shown.
By Condition S, $S_{\lambda_0}'$ 
does not intersect the $t$-axis after $a$ and $M$
does not intersect the $t-axis$ after $b$.
If $a=b$, then we represent, for $t>a$ but close to $a$,
$S_{\lambda_0}'$ by
$y=u(t)$ and $M$ by $y=v(t)$.
By Condition S, we know that
$\dot u, \dot v>0$ for $t>a$ and $t$ close to $a$.
Applying Lemma \ref{lem8}, $u\equiv v$ near $t=a$.
It is  now reduced to Case (i) and 
the symmetry of $M$ follows.

If $a<b$, then Case (i) cannot occur.
Thus, by the definition of $\lambda_0$, both $S_{\lambda_0}'$ and
$M$ must have a horizontal tangent
line at $Q$. Let, as shown in Fig. 9,
$[Q,R]$ be the largest flat segment on the top part of
$S_{\lambda_0}'$ and
$[Q, S]$ be the largest flat segment on the top part of
$M$.  Then applying
Lemma \ref{lem8} as above, $R$ and $S$
must be different.  Recall that Case (i) 
does not happen.
It is then clear that for all  $\lambda$  close 
to $\lambda_0$, $S_{\lambda}'$ still lie in $\overline U$, contradicting
to the definition of $\lambda_0$.
Theorem \ref{thm1} is established.

\vskip 5pt
\hfill $\Box$
\vskip 5pt

The remaining sections take up some further
one dimensional model problems.

\section{More results}

\begin{lem}
 Let $u$ be a $C^2$ function
 on an interval $(a,b)$, which is $C^1$ on the
 closure, and satisfies
 \begin{equation}
 \dot u(a)=0
 \label{4.1new}
 \end{equation}
 and 
 \begin{equation}
  \dot u>0 \ \mbox{on}\ (a, b).
   \label{inew}
    \end{equation}
    Let $v$ be a $C^2$ function
    on an interval $(\alpha, \beta)$ which
    is $C^1$ in the closure and satisfies
    \begin{equation}
    v(\alpha)\le u(b),\quad
    v(\beta)=u(a), \quad\mbox{and}
    \ \ v(\beta)<v(s)<v(\alpha),\ \forall\ \alpha< s< \beta,
    \label{4.3a}
    \end{equation}
    and
    \begin{equation}
    \dot v(\alpha)=0.
    \label{4.4new}
    \end{equation}
    See for example

***********

Fig. 10

***********

    Suppose that
    $$
    \mbox{whenever}\ u(t)=v(s)\ \mbox{for some}\
    \alpha<s<\beta \
    \mbox{we have}\ \ddot u(t)\le v''(s).
    $$
    Then
     $v$ is a reflection of $u$:
     $v(t)\equiv u(c-t)$ where $c=b+\alpha=a+\beta$.  In particular,
      $v(\alpha)=u(b)$.
      \label{lem2new}
      \end{lem}

Lemma \ref{lem2new} is equivalent to

\noindent{\bf Lemma \ref{lem2new}$'$.}\
{\it In the hypotheses of Lemma \ref{lem2new}, if we change
(\ref{4.1new}) and (\ref{inew}) to}
$$
\dot u(b)=0
$$
{\it and}
 $$
 \dot u<0 \ \mbox{on}\ (a, b),
$$
      {\it
      and change (\ref{4.3a}) to
     }
      $$
      v(\alpha)\le u(a),\quad
      v(\beta)=u(b), \quad\mbox{and}
      \ \ v(\beta)<v(s)<v(\alpha),\ \forall\ \alpha< s< \beta,
     $$
      {\it Then $v$ is simply a translate of $u$.
      }
      
      \bigskip
      
\noindent{\bf Proof of the equivalence of Lemma  \ref{lem2new}
and  Lemma \ref{lem2new}$'$.}\
Let $U(t)=u(-t)$ and $V(t)=v(t)$.

\vskip 5pt
\hfill $\Box$
\vskip 5pt

\noindent{\bf Proof of Lemma  \ref{lem2new}.}\
Reflect $v$ about $\alpha$: Set
$$
w(t):= v(2\alpha-t), \quad
2\alpha-\beta<t<\alpha.
$$
Shift $w$ far to the right and then slide it
to the left until the graph first touches that
of $u$.  If the touching  occurs at $(c, u(c))$ for some $
a<c<b$,
then, by Lemma \ref{lema1},
the shift of $w$ coincides with $u$ neat $c$, which in turn 
implies that they coincide everywhere and $v$ is a reflection of $u$ 
as desired.
Since $\dot w(\beta)=\dot v(\alpha)=0$ while $\dot u>0$ on $(a,b)$, 
there are only two possibilities: The above situation does not
occur but the touching occurs at 
$(a, 0)$ or at $(b, u(b))$. 
If the touching occurs at $(a, 0)$, then
we must have $\dot v(\beta)=0$ since
$\dot u(a)=0$.
By Theorem \ref{prop1}
and Lemma \ref{lema1}, the two graphs must be identical near the origin.
Impossible.
If the touching  occurs at  $(b, u(b))$, then we must have
$\dot u(b)=0$ since we know that $\dot v(\alpha)=0$.
Let 
$\overline w$ denote the shift, and we know, for some $\epsilon>0$, that
$$\overline w(b)=v(\alpha)=u(b),
\ \ \bar w<u(b)\ \mbox{on}\ (b-\epsilon, b).
$$
Turning the picture upside down, and
applying Theorem \ref{prop1}, we again get a contradiction.
More precisely, let
$$
U(t):=
\overline w(t)-\overline w(b-t),\quad
V(t)=u(b)-u(b-t), \qquad 0< t<\epsilon.
$$
Applying Theorem \ref{prop1} to $U$ and $V$ leads
to $U\equiv V$ near the origin, i.e. $u\equiv \overline w$ 
near $b$.  Impossible.
Lemma  \ref{lem2new} is established.

\vskip 5pt
\hfill $\Box$
\vskip 5pt

Lemma \ref{lem2new}  is also equivalent to

\noindent{\bf Lemma \ref{lem2new}$''$.}\
{\it 
Let $u$ be a $C^2$ function on an interval
$(a,b)$, $C^1$ in the closure, satisfying
}
$$
\dot u(b)=0\ \ \ 
\mbox{and}\ \ 
u(a)>u(t)>u(b)\ \mbox{for}\ 
a<t<b,
$$
{\it
and $v$ be a $C^2$ function on an interval
$(\alpha, \beta)$, $C^1$ in the closure,
satisfying
}
$$
\dot v(\beta)=0,
$$
$$
\dot v>0,
$$
{\it and }
$$
v(\alpha)\le u(b), \quad v(\beta)=u(a).
$$
{\it 
Finally, assume that
}
$$
\mbox{whenever}\ u(t)=v(s)\
\mbox{we have}\ \ddot u(t)\le v''(s).
$$
{\it 
then $v$ is simply a reflection of $u$
and $v(\alpha)=u(b)$.
}

\bigskip

\noindent{\bf Proof of the equivalence of Lemma \ref{lem2new}
and  Lemma \ref{lem2new}$''$.}\
Let $V(t)=-u(-t)$ and
$U(t)=-v(-t)$.

\vskip 5pt
\hfill $\Box$
\vskip 5pt

\begin{thm} Let $u$ be a positive $C^2$ function on $(0,b)$ with
$u\in C^1([0,b])$,  satisfying
\begin{equation}
u(0)=\dot u(0)=u(b)=0.
\label{2.3a}
\end{equation}
Let $a$ be the first point where $u$ achieves its maximum
and assume that
\begin{equation}
\dot u>0\quad \mbox{on}\ (0,a).
\label{2.4a}
\end{equation}
Assume furthermore (main condition) that
\begin{equation}
\mbox{whenever}\ u(t)=u(s)\ \mbox{for}\
t<s\ \mbox{also}\ \ddot u(t)\le u''(s).
\label{2.5a}
\end{equation}
$\underline{\mbox{\it Conclusion}}$: $u$ is symmetric
about $b/2$ and 
\begin{equation}
u\equiv u(a)\ \ \mbox{on}\ [a, b-a].
\label{2.6a}
\end{equation}
\label{prop2}
\end{thm}

Note that we do not assume that $\dot u(b)=0$.

\begin{rem} The condition $\dot u(0)=0$ cannot be dropped.  Indeed,
we could consider a positive symmetric function on some
interval $(0,b)$ satisfying (\ref{2.3a}) except
that $\dot u(0)>0$.  Then, near $b$
we could change $u$ slightly by increasing its 
second derivative there, in such a way that the new function, when extended
would still vanish at some point $\bar b>b$.  The resulting function on
$(0,\bar b)$ would satisfy (\ref{2.4a}) but would not be symmetric.
\label{rem2.1}
\end{rem}

Here is an example showing that if condition (\ref{2.4a})
is weakened to $\dot u\ge 0$ on $(0, a)$, then
$u$ need not be symmetric.  Here
$u$ on $(4, 5)$ is the reflection
of $u$ on $(0, 1)$.

***********

Fig. 11

***********

\noindent{\bf Proof of Theorem  \ref{prop2}.}\
Let $b_1$ be the last value of $t$ where $u$ assumes
its maximum. By Lemma \ref{lem2new} ,
$b_1=b-a$ and
$$
u(t)=u(b-t)\quad \mbox{for}\ 0\le t\le a.
$$
Now we prove that $u$ is constant on
$[a, b_1]$.  If not, we can find
$[\alpha, \beta]\subset [a,  b_1]$ such that
$$
u(a)\ge u(\alpha)>u(t)>u(\beta)=
\min_{ [a, b_1] }u>0,
\qquad \mbox{for}\ \alpha< t<\beta.
$$

Since $u(t)=u(b-t)$ on $(0, a)$,
it follows from the main condition, that if $t\le a<s$ and
$u(t)=u(s)$ then
$\ddot u(t)=u''(s)$.
Thus for $s$ on $(\alpha, \beta)$ we can
find unique $t(s)$ on $(0, a)$
such that
$$
u(t(s))=u(s).
$$
Hence
$$
\dot u(t(s))t'=u'(s)
$$
and
$$
\left( \dot u(t(s))^2 -u'(s)^2\right)'
=2 u'(s)\left( \ddot u(t(s))-u''(s)\right)
=0.
$$
Hence
$$
\dot u(t(s))^2 =u'(s)^2.
$$
This is impossible, since $u'(\beta)=0$, while
$\dot u(t(\beta))>0$.
Theorem \ref{prop2} is proved.

\vskip 5pt
\hfill $\Box$
\vskip 5pt

\section{General second order operators}

In this section we extend various results
 to nonlinear second order ordinary differential equations.
We consider 
\begin{equation}
K\in C^0(\Bbb R^3),\ K(s,p,q)\ \mbox{is } C^1 \mbox{ in}\ (p,q),\ 
 \mbox{and}\
K_q(s,p,q)>0,\ \forall\ (s,p,q)\in \Bbb R^3,
\label{K1}
\end{equation}
and we study nonlinear second order differential operator
$K(u, \dot u, \ddot u)$.  It is elliptic because of (\ref{K1}).

The first result is an extension of Lemma \ref{lema1}.
\begin{lem} Let $K$ satisfy (\ref{K1}), and
let $u\ge v$ be two $C^2$ functions on
$(0,b)$
satisfying
\begin{equation}
\max\{\dot u, \dot v\}>0\ \mbox{on}\ (0,b),
\label{aa0z}
\end{equation}
and
\begin{equation}
\mbox{if}\ u(t)=v(s)\ \mbox{for}
\ 0<t<s<b\
\mbox{there}\ 
K(u(t), \dot u(t), \ddot u(t))\le K(v(s), v'(s), v''(s)).
\label{2.2z}
\end{equation}
Then either
$$
u>v\  \mbox{on}\ (0,b),
$$
or
\begin{equation}
u\equiv v\ \ \mbox{on}\ (0, b).
\label{aa1z}
\end{equation}
\label{lema1z}
\end{lem}

\noindent{\bf Proof.}\
Suppose $u(c)=v(c)$ for some $0<c<b$. Then, by (\ref{aa0z}),
$$
\dot u(c)=\dot v(c)>0.
$$
By  the
implicit function theorem, for $s$ close to $c$, there is a $C^2$ function $t(s)$ such that
\begin{equation}
u(t(s))=v(s).
\label{2.4z}
\end{equation}
Set, for $s$ close to $c$,
$$
g(s)=s-t(s),
$$
so $g\ge 0$.  Differentiating
(\ref{2.4z}) we find
$$
\dot u(t(s)) t'=v'(s),
\quad \ddot u t'^2 +\dot u t''=v''.
$$
Thus, for some functions $c_1(s)$ and $c_2(s)$,
we have
$$
\dot u(t(s))-v'(s)=c_1(s)g'(s),\quad
\ddot u(t(s))-v''(s)
=\dot u g''+c_2(s)g'.
$$
Using (\ref{K1}) and the above, 
we obtain, via the mean value theorem,
\begin{equation}
0\ge  K(u, \dot u, \ddot u)-K(v, v', v'')
=a(s) \dot u g''+ c(s)g',\ \ \mbox{with}\ a(s)>0,
\label{2.5z}
\end{equation}

Now if $u(s)=v(s)$ for some $s>0$
then $g$ vanishes there.  By the strong maximum
principle applied to (\ref{2.5})
 implies that $v\equiv u$ in a neighborhood of $c$. 
(\ref{aa1z}) follows immediately.
Lemma \ref{lema1z} is established.

\vskip 5pt
\hfill $\Box$
\vskip 5pt

The second result is an extension of Lemma \ref{maximum}.
\begin{lem} Let $K$ satisfy (\ref{K1}), and let
$u, v\in C^2((0,b))\cap C^0([0, b])$
satisfy
$$
v(0)\le u(0),\ v(b)\le u(b),\
\mbox{and}\ v< u(b)\ \mbox{on}\ (a,b),
$$
$$
\mbox{either}\ \dot u>0\  \mbox{on}\ (0, b)\
\mbox{or}\ \dot v>0 \ \mbox{whenever}\ u(0)<v<u(b),
$$
and
(\ref{2.2z}).
Then
$$
u\ge v\ \mbox{on}\ [0,b].
$$
\label{maximumz}
\end{lem}

\noindent{\bf Proof.}\  The proof is essentially the same as that of
Lemma \ref{maximum}.  The only difference is to use
Lemma \ref{lema1z} instead of Lemma \ref{lema1}.

\vskip 5pt
\hfill $\Box$
\vskip 5pt

The third result is an extension of Lemma \ref{lem1}.   

\begin{lem} Let 
 $u$ and $w$ be positive $C^{1,1}$ functions
in $(0, c)$, belonging to $C^1([0, c))$, and both satisfying
\begin{equation}
\frac{d}{dt} \left( K(\dot u) \right)=f(u),
\label{3.3z}
\end{equation}
and
$$
u(0)=w(0)=0,
\quad
\dot u(0)=\dot w(0).$$
Here
\begin{equation}
K\in C^1(\Bbb R), \ \mbox{and}\
K'>0\ \mbox{in}\ \Bbb R.
\label{K2}
\end{equation}
Assume in addition that $K'$ is even, and assume
 that $f(\rho)$ is locally Lipschitz for
$\rho>0$ (not necessarily on $\rho\ge 0$).
Then
$$
u\equiv w.
$$
\label{lem1y}
\end{lem}

\noindent{\bf Proof.}\  The proof is similar to that of Lemma \ref{lem1}.
Multiplying (\ref{3.3z}) by $\dot u$ we find
that
\begin{equation}
\frac{d}{dt} G(\dot u)= \frac{d}{dt} F(u)
\label{bb2}
\end{equation}
where
\begin{equation}
G(p)=\int_0^p \rho K'(\rho)d\rho,\ \ \mbox{and}\
F\ \mbox{is such that}\ F_u=f(u).
\label{bb3}
\end{equation}
 Integrating this from $t_0$ to $t$,
    $t_0>0$, we find
    \begin{equation}
    G(\dot u(t))-G(\dot u(t_0))=
    F(u(t))-F(u(t_0)).
    \label{bb1}
    \end{equation}
Letting $t_0\to 0$ we see that
$F(u)$ has a limit at $u=0$, which we may take
to be zero.
Thus
\begin{equation}
G(\dot u)=F(u).
\label{3.5z}
\end{equation}
As in the proof of  Lemma \ref{lem1}, we have
$\dot u>0$ on $(0, \frac c2)$.
Since $K'>0$, we see that
\begin{equation}
G(p)>0, \ G'(p)>0,\qquad \mbox{for}\ p>0.
\label{posi}
\end{equation}
Thus
$$
F(u(t))>0,\qquad \mbox{for}\ 0<t<\frac c2.
$$
It follows that
$$
\frac {\dot u}
{ G^{-1}(F(u(t)) }=1,\qquad \mbox{on}\ (0, \frac c2).
$$
If, on $\rho>0$, $H(\rho)$ is such that
$$
H_\rho=\frac 1{   G^{-1}(F(\rho))  },
$$
we find  that
$$
\frac {d}{dt} H(u)=1.
$$
Integrating from $t_0>0$ to $t>t_0$ we obtain
$$
H(u(t))-H(u(t_0))=t-t_0.
$$
Letting again $t_0\to 0$
we see that $H$ has a limit at $\rho=0$, which we may take to be
zero.  Thus
$$
H(u(t))=t.
$$
So $u(t)$ is uniquely determined on $(0, c/2)$,
since $H_\rho>0$ for $\rho>0$.
Then, by the local Lipschitz property of $f(\rho)$ on
$\rho>0$, it follows that $u$ is unique on
$(0, c)$.  Lemma
\ref{lem1y} is proved.

\vskip 5pt
\hfill $\Box$
\vskip 5pt

The fourth result is an extension of Lemma \ref{lem1prime}.
Here we do not assume that $K'$ is even.
\begin{lem}
Let $K$ satisfy (\ref{K2}), and
let $u$ and $w$ be positive $C^{1,1}$ functions
on $(0,c)$, belonging to $C^1( [0, c))$, and satisfying,
for some $f\in L^\infty(0,\infty)$,
$$
\frac{d}{dt} \left( K(\dot u) \right)=f(u), \ \ \
\frac{d}{dt} \left( K(\dot w) \right)=f(w),\qquad \mbox{in}\ (0,c),
$$
and
$$
u(0)=w(0)=0,
\ \ \
\dot u(0)=\dot w(0),\ \ \
\dot u>0\ \mbox{on}\ (0,c).
$$
Then
$$
u\equiv w\ \mbox{on}\ (0,c).
$$
\label{lem5-4}
\end{lem}

\noindent{\bf Proof.}\
Follow the proof of Lemma  \ref{lem1y} until (\ref{3.5z}).
In a similar way we also have
\begin{equation}
G(\dot w(t))=F(w(t))\qquad \mbox{on}\ (0,c).
\label{abc1new}
\end{equation}
Let $0<b<c$ be any number satisfying
(\ref{number}). We  
only need to prove (\ref{A4-1}).

By (\ref{3.5z}), (\ref{abc1new})  and the fact that $\dot u>0$ on $(0, c)$,
we know that
$$
G(\dot w(t))=F(w(t))>0,
\qquad\forall\ 0<t<b.
$$
Since $w(0)=0$ and $w>0$
on $(0, b)$, we have $\dot w>0$ on $(0, b)$.
Proceeding as in the proof of Lemma  \ref{lem1y}, we arrive at
$$
H(u(t))=H(w(t))=t,\qquad
0<t<b.
$$
But $H_\rho(\rho)>0$ for $\rho>0$  and $H(0)=0$, we obtain
(\ref{A4-1}).  Lemma \ref{lem5-4} is established.

\vskip 5pt
\hfill $\Box$
\vskip 5pt

The fifth result is an extension of Lemma \ref{lem8}. 
\begin{lem} Let $K$ satisfy (\ref{K2}), and let, 
 for some $b>0$,  $u, v\in C^2( (0, b))\cap C^1([0,b])$
 satisfy
 $$
 u(t)\ge  v(t)>0,\qquad  \mbox{for}\
  0<t\le b,
  $$
  $$
  \mbox{either}\ 
 \dot u(t)>0,   \dot v(t)\ge 0\ \mbox{for}\
  0<t\le b,
  \ \mbox{or}\
  \dot u(t)\ge 0,   \dot v(t)>0,\ \mbox{for}\
 0<t\le b,
 $$
 and
 $$
 u(0)=\dot u(0)=0.
 $$
 Assume
 \begin{equation}
 \mbox{whenever}\ u(t)=v(s)\ \mbox{for}\
 0<t<s<b\ \mbox{we have}\
 \frac{d}{dt}
 \left(K(\dot u(t))\right)\le
  \frac{d}{ds}
   \left(K(\dot v(s))\right).
   \label{3.1z}
   \end{equation}
 Then
 \begin{equation}
 u\equiv v\ \ \mbox{on}\ [0, b].
 \label{3.2newz}
 \end{equation}
 \label{lem8z}
 \end{lem}

\noindent{\bf Proof.}\  
The proof is similar to that of Lemma \ref{lem8}. 
 We will only prove it under
 ``$ \dot u(t)>0, \dot v(t)\ge 0$ for
 $0<t\le b$'', since the changes needed when
 assuming instead
 ``$\dot u(t)\ge 0,  \dot v(t)>0$ for
 $0<t\le b$'' are similar to those
 in the proof of Theorem \ref{prop1}.
We start as in the proof of Theorem 
\ref{prop1}.  $u$ satisfies (\ref{3.3z})
for some unknown  continuous  function $f$  on $[0,  u(b)]$.

Condition (\ref{3.1z}) means
\begin{equation}
     K'(\dot v(s)))\ddot v \ge f(v).
\label{3.4z}
\end{equation}
Multiply (\ref{3.3z}) by $\dot u$; we find (\ref{bb2})
with $G$ given by (\ref{bb3}).
(\ref{3.5z}) still holds, so does (\ref{posi}).

Multiplying  (\ref{3.4z}) by $\dot v\ge 0$ we obtain
$$
\frac d{dt}G(\dot v)\ge
 \frac d{dt} F(v) \ \
\mbox{for}\ 0<t_0<t.
$$

Since (\ref{posi}) still holds, and since $\dot v\ge 0$, we have
$$
G(\dot v(t))\ge 0,\ \mbox{for}\ 0\le t\le b.
$$
Thus, also in view of 
 our choice of setting $F(0)=0$,
 we find
by integrating
\begin{equation}
G(\dot v)\ge F(v).
\label{3.6z}
\end{equation}
Because of the the second inequality in
(\ref{posi}),  (\ref{3.5z}) and (\ref{3.6z}) imply that
(\ref{3.7}), and the rest of the proof follows exactly as
the arguments after (\ref{3.7}) in the proof of Lemma \ref{lem8}.
Lemma \ref{lem8z} is established.

\vskip 5pt
\hfill $\Box$
\vskip 5pt

The sixth result is closely related to Lemma \ref{lem2new}.
\begin{lem}
Let $K$ satisfy (\ref{K2}), and let
 $u$ be a $C^2$ function
  on an interval $(a,b)$, which is $C^1$ on the
   closure, and satisfies (\ref{4.1new}) and
   (\ref{inew}).
     Let $v$ be a $C^2$ function
         on an interval $(\alpha, \beta)$ which
	     is $C^1$ in the closure and satisfies
	     (\ref{4.3a}), (\ref{4.4new}) and
	     \begin{equation}
	     \dot v\le 0,\qquad \mbox{on}\ (\alpha, \beta).
	     \label{extra}
	     \end{equation}
	      Suppose that
	          $$
		      \mbox{whenever}\ u(t)=v(s)\ \mbox{for some}\
		          \alpha<s<\beta \
			      \mbox{we have}\ 
			      \frac{d}{dt}\left(
			      K(\dot u(t))\right)
			      \le  \frac{d}{ds}\left(
			                                    K(\dot 
v(s)))\right).
$$
  Then
       $v$ is a reflection of $u$:
            $v(t)\equiv u(c-t)$ where $c=b+\alpha=a+\beta$.  In particular,
	          $v(\alpha)=u(b)$.
		        \label{lem5-6}
			\end{lem}

\begin{rem}
There are analogues of  Lemma \ref{lem2new}$'$,  \ref{lem2new}$''$,
which we call  Lemma \ref{lem5-6}$'$,  \ref{lem5-6}$''$, and which follow from
 Lemma \ref{lem5-6} as do the others from  Lemma \ref{lem2new}.
 \end{rem}

			\noindent{\bf Proof of Lemma \ref{lem5-6}.}\
Follow the proof of  Lemma  \ref{lem2new} and
make the following changes:
Change ``by Lemma \ref{lema1}'' to ``by Lemma 
\ref{lema1z}''; change
`` by Theorem \ref{prop1}
and Lemma \ref{lema1}'' to ``by 
 Lemma \ref{lem8z}'';
 change ``applying  Theorem \ref{prop1}'' to ``applying  Lemma \ref{lem8z}''.

\vskip 5pt
\hfill $\Box$
\vskip 5pt

\begin{question}
Do the conclusions of Lemma \ref{lem1y},
Lemma \ref{lem5-4}, Lemma \ref{lem8z} and Lemma \ref{lem5-6}
still hold if we replace
$\displaystyle{
 \frac{d}{dt}
    \left(K(\dot u(t))\right)\le
           \frac{d}{ds}
	                 \left(K(\dot v(s))\right)}$
			                      in (\ref{3.1z}) by
                                  $K(u(t), \dot u(t), \ddot u(t))\le
                                   K(v(s), v'(s), v''(s))$ for some $K$
													                                               satisfying
																		                                                        (\ref{K1}), or even for those $K(s,p,q)$ which are
																									                                                        independent of $s$?
																																                                                               \end{question}

\section{Counter examples}

We will give an example showing that
the answer to Question \ref{question2} on page ? is no.

First, we present a function $u$ on $(0, 4)$.

***********

Fig. 12

***********

Here, $u$ on $(3,4)$ is the reflection
of its values on $(0,1)$.  $u$ will satisfy the
curvature condition:
\begin{equation}
\mbox{for any  } t<s\
\mbox{such that}\ u(t)=u(s),
\ 
\frac{ \ddot u(t) }
{  (1+\dot u(t)^2 )^{ \frac 32}  }
\le 
\frac{  u''(s) }
{  (1+u'(s)^2)^{ \frac 32}  }
\ \mbox{hold}.
\label{Q64}
\end{equation}

$u$ will be taken to be symmetric on
$(1,2)$ about $\frac 32$, and symmetric
on $(2,3)$ about $\frac 52$.  We will then require
(\ref{Q64}) only for $1\le t\le \frac 32$ and
$2\le s\le \frac 52$.

For convenience, after subtracting a constant, and shifting,
we may describe the two bumps by two functions $u$ and $v$
(we still call the first $u$) on
$(0, 1)$ given by 
\begin{equation}
u=\epsilon^6 t^3 (1-t)^3,\quad
v=\epsilon^3 t^3(1-t)^3
\label{4.14}
\end{equation}
with $\epsilon$ very small.

\noindent $\underline{  \mbox{Claim}}$.\
Whenever $u(t)=v(s)$ for $t<\frac 12$, there
$\displaystyle{
\frac{\ddot u(t) } { (1+\dot u(t)^2)^{\frac 32} }\le
\frac{ v''(s) }  { (1+v'(s)^2)^{\frac 32} }
}$.

\noindent $\underline{  \mbox{Proof}}$.\
$u(t)=v(s)$ means
\begin{equation}
\epsilon t(1-t)=s(1-s).
\label{4.15}
\end{equation}

  We have
only to check for the function,
$u$ and $v$ in (\ref{4.14}) that
\begin{equation}
\mbox{if}\ u(t)=v(s)\ \mbox{for}
\ 0<t, s<\frac 12\ \mbox{there}\
\frac{\ddot u(t) } { (1+\dot u(t)^2)^{\frac 32} }\le  
\frac{ v''(s) }  { (1+v'(s)^2)^{\frac 32} }.
\label{4.16}
\end{equation}

\noindent $\underline{  \mbox{Proof of (\ref{4.16})}}$.\
From (\ref{4.15}) we see that for $\epsilon$ small,
\begin{equation}
s=\epsilon (t-t^2)+
O\left(  \epsilon^2(t-t^2)^2\right).
\label{4.17}
\end{equation}
Now
\begin{equation}
\dot u = \epsilon^6 3 (1-2t)(t-t^2)^2,
\label{4.18}
\end{equation}
\begin{equation}
\ddot u=6\epsilon^6
\left[ - (t-t^2)^2 + (1-2t)^2 (t-t^2)\right]=
6\epsilon^6    (t-t^2) ( 1-5t +5t^2).
\label{4.19}
\end{equation}
Thus
$$
\frac {\ddot u}{  (1+\dot u^2)^{\frac 32} }
\le C\epsilon^6 (t-t^2).
$$
At the same time
$$
\frac {v''}{  (1+(v')^2)^{\frac 32} }
=\epsilon^3 \left[ 6s+O(s^2)\right]\ge 5\epsilon^4(t-t^2),
$$
by (\ref{4.17}).
It follows that (\ref{4.16}) holds for
$\epsilon$ small.

\vskip 5pt
\hfill $\Box$
\vskip 5pt

Finally we obtain an example of a closed nonconvex curve
$M$ satisfying $H(A)\le H(B)$
if $A, B$ lie on $M$ and $A_1\le B_1$,
but with $H(A)$ not equal to $H(B)$.

Namely take the curve $u$ above and round it off on the bottom in
a symmetric way.

***********

Fig. 13

***********

\section{Open problems in higher dimension}

The problems are related to Theorem 
\ref{prop1} and to
Lemma \ref{lem1prime}, \ref{lem1}.
For convenience, we will denote the 
independent variables by $(t, y)$, $t$ nonnegative,
$y$ in $\Bbb R^{n-1}$.
The functions we consider are defined in the closure
of the half ball
$$
B^+:=\{(t, y)\ |\ t^2+|y|^2<R^2, \ t>0\}.
$$

The first question is related to Theorem \ref{prop1}.

\begin{problem}
Suppose $u\ge v>0$ and $\displaystyle{
u_t:=\frac{\partial u }{ \partial t}
>0
}
$ in $B^+$, and $u$ and $v$ are 
$C^2$ in the closure of $B^+$.  Suppose  
\begin{equation}
u=v=u_t
=v_t=0\ 
\mbox{on}\ \{(0, y)\ |\ |y|<R\}.
\label{70}
\end{equation}
Assume the main condition: 
$$
\mbox{whenever}\ 
u(t,y)=v(s,y)\ \mbox{for}\ t<s,
\mbox{then}\ \Delta u(t,y)\le \Delta v(t,y)\ \mbox{holds}.
$$
Question:  Is $u\equiv v$?
\label{problem7.1}
\end{problem}

\begin{problem}
Let $u$ and $v$ be $C^\infty$ in the closure 
of $B^+$ and positive in $B^+$, and satisfy
$
\Delta u=f(y, u) 
$, $\Delta v=f(y, v)$ and satisfying
(\ref{70}).
Concerning $f$, we  assume that it is
continuous in $u\ge 0$, smooth
in $y$ there, and it is smooth in $(y,u)$ where $u>0$.

Question: Is $u\equiv v$? 
\label{problem7.2}
\end{problem}

In Part II we prove 
the answer is yes, but under additional assumptions, that $u$ and $v$
vanish of finite order in $t$ at the
origin and that $u\ge v$.

\begin{problem}
Is the answer to Problem \ref{problem7.2} yes if we add the hypothesis that
$u_t>0$ when $t>0$?
\label{problem7.3}
\end{problem}


\begin{thebibliography}{99}
\bibitem{A} A.D. Alexandrov, Uniqueness theorems for surfaces in the large,
 Vestnik,
Leningrad Univ. 13 (1958), 5-8. 
 \bibitem{L} Y.Y. Li,
 Group invariant convex hypersurfaces with prescribed
 Gauss-Kronecker curvature, Contemporary Mathematics,
 AMS, 205 (1997), 203-218.
\bibitem{LN} Y.Y. Li and L. Nirenberg,
 A geometric problem and the Hopf Lemma. II, in preparation.
 \bibitem{M} J. Mawhin,
{\it  Points fixes, points critiques et problèmes aux limites},
 S\'eminaire de Mathématiques Sup\'erieures
 92,
 Presses de l$'$Universit\'e de Montr\'eal, Montreal, QC, 1985.
\bibitem{Ros} A. Ros, Compact hypersurfaces with 
constant higher order mean curvatures,
Rev. Mat. Iberoamericana 3 (1987),  447-453.
\bibitem{W} H.C. Wente, Counterexample to a conjecture of 
H. Hopf,
 Pacific J. Math. 121 (1986), 193-243. 
 \end{thebibliography}
\end{document}